\documentclass[12pt,a4paper,reqno]{amsart}
\usepackage{amsmath}
\usepackage{amsfonts}
\usepackage{amssymb}
\numberwithin{equation}{section}

\theoremstyle{plain}

\newtheorem{theorem}[subsection]{Theorem}
\newtheorem{proposition}[subsection]{Proposition}
\newtheorem{lemma}[subsection]{Lemma}

\newtheorem{conjecture}[subsection]{Conjecture}

\theoremstyle{definition}

\newtheorem{definition}[subsection]{Definition}
\newtheorem{problem}[subsection]{Problem}
\newtheorem{example}[subsection]{Example}
\newtheorem{question}[subsection]{Question}

\renewcommand{\leq}{\leqslant}
\renewcommand{\geq}{\geqslant}

\newsavebox{\proofbox}
\savebox{\proofbox}{\begin{picture}(7,7)%
  \put(0,0){\framebox(7,7){}}\end{picture}}
\def\boxeq{\tag*{\usebox{\proofbox}}}

\def\proof{\noindent\textit{Proof. }}
\def\endproof{\hfill{\usebox{\proofbox}}}

\def\Agh{A_H^{+g}}
\def\ind{\mbox{ind}}

\begin{document}

\title{Finite Field Models in Arithmetic Combinatorics
}

\author{Ben Green}
\address{Trinity College\\
     Cambridge CB2 1TQ\\
     England
}
\email{bjg23@hermes.cam.ac.uk}

\thanks{The author is a Fellow of Trinity College, Cambridge.}

\begin{abstract}
The study of many problems in additive combinatorics, such as Szemer\'edi's theorem on arithmetic progressions, is made easier by first studying models for the problem in $\mathbb{F}_p^n$, for some fixed small prime $p$. We give a number of examples of finite field models of this type, which allows us to introduce some of the central ideas in additive combinatorics relatively cleanly. We also give an indication of how the intuition gained from the study of finite field models can be helpful for addressing the original questions.

\end{abstract}

\maketitle

\section{Introduction}\label{sec1}
\noindent This article is concerned with a variety of problems in additive and combinatorial number theory. The following two examples will convey the general flavour:
\begin{problem}[3-term APs]\label{prob1} What is $r_3(N)$, the cardinality of the largest subset of $\{1,\dots,N\}$ containing no three distinct elements $x,x+d,x+2d$ in arithmetic progression?\end{problem}
\begin{problem}[Sets with small doubling]\label{prob2} If $A \subseteq \mathbb{Z}$, write $A + A$ for the set of all sums $a + a'$, $a,a' \in A$. What can be said about the structure of $A$ if $A$ is nearly closed under addition in the sense that $|A + A| \leq K|A|$?\end{problem}
\noindent What, then, is the ``general flavour''? Of course, both of these problems are of an additive combinatorial flavour. Furthermore, they may both be asked in a general abelian group. Regarding Problem \ref{prob1}, we may define the quantity $r_3(G)$ for any finite abelian group $G$. And Problem \ref{prob2} makes sense in any abelian group.\vspace{11pt}

\noindent The ability to generalise to an arbitrary $G$ will be a common feature of many of the questions we discuss. An important observation is that not all abelian groups were created equal. It turns out that both Problems \ref{prob1} and \ref{prob2} are both considerably easier in groups other than those in which they were originally asked ($\mathbb{Z}/N\mathbb{Z}$ for Problem \ref{prob1}\footnote{In many questions, the difference between $\mathbb{Z}/N\mathbb{Z}$ and $\{1,\dots,N\}$ is purely technical.} and $\mathbb{Z}$ for Problem \ref{prob2}). 
Indeed, Meshulam \cite{meshulam} observed that Problem \ref{prob1} is naturally addressed in $\mathbb{F}_3^n$, whereas Ruzsa \cite{ruzsa-freimangroups} saw that Problem \ref{prob2} is particularly pleasant in $\mathbb{F}_2^{\infty}$.  Here, $\mathbb{F}_p$ denotes the finite field with $p$ elements, and $\mathbb{F}_p^{\infty}$ is our notation for a vector space of countable dimension over $\mathbb{F}_p$.\vspace{11pt}

\noindent Roughly speaking, the reason that finite field models are nice to work with is that one has the tools of linear algebra, including such notions as \textit{subspace} and \textit{linear independence}, which are unavailable in general abelian groups.\vspace{11pt}

\noindent Historically, questions such as Problems \ref{prob1} and \ref{prob2} were investigated in their original settings, and it was observed only later that analogous arguments worked in the finite field setting and in fact looked rather simpler. More recently, there has been a trend in the opposite direction. This has been fuelled by an idea of Bourgain \cite{bourgain2} which, suitably interpreted, can be viewed as a way of converting arguments in the finite field setting to arguments which work for an arbitrary group $G$ by using a kind of ``approximate linear algebra''. The author \cite{green-reg} produced a result about sets of integers with few solutions to $x + y = z$ which would have been very difficult to attain without first considering a finite field model, and more work of this sort is in progress. It is an interesting feature of many problems that progress for the groups $G$ which are ``of interest'', such as $\mathbb{Z}/N\mathbb{Z}$, is scarcely simpler than for general abelian $G$.\vspace{11pt}

\noindent The format of this article is as follows. After setting up a little notation and a few definitions, we will discuss a number of finite field problems of ``Szemer\'edi type'', that is to say along the lines of Problem \ref{prob1}. We will strive for a uniform treatment of three such problems: 3-term APs (\S \ref{sec4}), right-angled triangles (\S \ref{sec5}) and 4-term APs (\S \ref{sec6}). We will discuss a fourth problem in \S \ref{sec7}, which concerns solutions to $x + y = z$ and is in a somewhat similar spirit.\vspace{11pt}

\noindent After these four sections we will, in \S \ref{sec8}, sketch an argument of Bourgain, which is currently being developed by the author and others, including T.Tao and I. Shkredov, into a machine for converting arguments in the finite field setting into arguments that work in any finite abelian group $G$. This is often of some interest when $G = \mathbb{Z}/N\mathbb{Z}$, because in that case it is often possible to infer results concerning the integers. \vspace{11pt}

\noindent After that there follow three further sections of a somewhat miscellaneous nature dealing with finite field analogues of problems in additive number theory.\vspace{11pt}

\noindent Since this is a survey article we have not gone into a great deal of technical detail. There are, however, two areas we discuss which are not well covered in the literature. Thus on the author's webpage one may find two supplementary documents \cite{sup-doc-1, sup-doc-2}. The first of these gives details of the finite field version of Shkredov's argument, which is outlined in \S \ref{sec5}. The second supplies proofs for the result of Ruzsa discussed in \S \ref{sec10}.\vspace{11pt}

\noindent Our scope in this article is a little limited, in that our main interest is in additive combinatorial problems which can be usefully studied in $\mathbb{F}_p^n$ for fixed $p$, regarding $n$ as a variable parameter. Secondly, I have unashamedly prioritised areas in which I have personally worked. There are most assuredly other areas of mathematics where finite field models have proved invaluable, such as the study of the Kakeya and restriction phenomena. We do not touch upon these matters here, referring the reader instead to the article \cite{mockenhaupt-tao} as well as in the surveys \cite{katz-tao-survey, tao-survey1, tao-survey2}.

\section{Notation and Basic Definitions}\label{sec2}
\noindent Let $p$ be a prime ($p$ will be either 2,3 or 5). Write $\mathbb{F}_p$ for the finite field with $p$ elements, which may be identified with $\mathbb{Z}/p\mathbb{Z}$, and for an integer $n \geq 1$ write $\mathbb{F}_p^n$ for a vector space of dimension $n$ over $\mathbb{F}_p$. This will be understood to have been given to us with a fixed basis $(e_1,\dots,e_n)$, relative to which we will occasionally write a given $x \in \mathbb{F}_p^n$ as a coordinate vector $(x_1,\dots,x_n)$. We will always write $N = p^n$ for the cardinality of the space $\mathbb{F}_p^n$.\vspace{11pt}

\noindent Once we have a basis the \textit{Fourier transform} of a function $f : \mathbb{F}_p^n \rightarrow \mathbb{C}$ can be written down in a concrete form. A complete set of characters $\gamma : \mathbb{F}_p^n \rightarrow S^1$ is given by the maps $\gamma_{\xi}$, defined by
\[ \gamma_{\xi}(x) = \gamma_{\xi_1,\dots,\xi_n}(x) = \omega^{\xi_1 x_1 + \dots + \xi_n x_n} = \omega^{\xi^{T}x},\] where $\xi \in \mathbb{F}_p^n$ and $\omega = e^{2\pi i/p}$. Thus, for any $\xi \in \mathbb{F}_p^n$, we define
\[ \widehat{f}(\xi) := \sum_x f(x) \gamma_{\xi}(x) = \sum_x f(x) \omega^{\xi^Tx}.\]
We may also write this as $f^{\wedge}(\xi)$ on occasion.
The basic facts concerning the Fourier transform are summarised in the following lemma.
\begin{lemma}[The Fourier Transform]\label{basic-fourier-properties}
Let $f,g : \mathbb{F}_p^n \rightarrow \mathbb{C}$ be two functions. Then
\begin{enumerate}
\item $\widehat{f}(0) = \sum_{x} f(x)$;
\item \textup{(Plancherel)} $\sum_x f(x)\overline{g(x)} = N^{-1} \sum_{\xi} \widehat{f}(\xi)\overline{\widehat{g}(\xi)}$;
\item \textup{(Inversion)} $f(x) = N^{-1} \sum_{\xi} \widehat{f}(\xi) \omega^{-\xi^T x}$;
\item \textup{(Convolution)} Write $(f \ast g)(x) = \sum_y f(y)g(x - y)$. Then $(f \ast g)^{\wedge}(\xi) = \widehat{f}(\xi)\widehat{g}(\xi)$.
\end{enumerate}
\end{lemma}
\noindent Very often, we will be concerned with functions $f$ which are the characteristic functions of sets $A \subseteq \mathbb{F}_p^n$. It is very convenient to abuse notation and write $A(x)$ for such a function. Thus $A(x) = 1$ if $x \in A$, and $A(x) = 0$ otherwise. This notation is by now reasonably widespread in the literature, as are alternative notations such as $\chi_A$ or $\mathbf{1}_A$.\vspace{11pt}

\noindent It will be very convenient to use the language of conditional expectation. Suppose that $x$ is a variable or set of variables, and that $f$ is a real-valued function of $x$. Then we write
\[ \mathbb{E}(f(x) | x \in B) := |B|^{-1}\sum_{x \in B} f(x)\]
for the average of $f(x)$ over all $x \in B$.\vspace{11pt}

\noindent Let us conclude with some notation concerning sumsets. If $G$ is an abelian group and if $A,B \subseteq G$ then we write $A + B = \{a + b : a \in A, b \in B\}$. For any positive integers $s,t$ we write $sA + tB$ for the set of all sums $a_1 + \dots + a_s + b_1 + \dots + b_t$,and $sA - tB$ for the set of all sums $a_1 + \dots + a_s - b_1 - \dots - b_t$.

\section{Uniformity}\label{sec3}

\noindent A notion which will feature repeatedly in this article is that of \textit{uniformity}, also referred to in various related guises as regularity, pseudorandomness or quasirandomness. 

\begin{definition}
Let $A \subseteq \mathbb{F}_p^n$ be a set, and let $\eta  \in (0,1)$ be a parameter. We will say that $A$ is $\eta$-uniform if 
\[ \sup_{\xi \neq 0} |\widehat{A}(\xi)| \leq \eta N.\]
\end{definition}
\noindent Observe that if $A$ is $\eta$-uniform then it is also $\eta'$-uniform for all $\eta' \geq \eta$.\vspace{11pt}

\noindent The basic philosophy behind this definition is as follows. A truly random set $A$ (generated, say, by including each $x \in \mathbb{F}_2^n$ in $A$ independently at random with probability $1/2$) will be $\eta$-uniform with very high probability. In fact, using a large deviation estimate such as Chernoff's bound (see \cite{alon-spencer} for example) one can show that this is true even for $\eta = N^{-1/2 + \epsilon}$. A truly random set will have many other properties almost surely. Remarkably, many of these are consequences of $A$ being $\eta$-uniform. This phenomenon was investigated in the context of graphs by Thomason \cite{thomason1,thomason2} and by Chung, Graham and Wilson \cite{chung-graham-wilson}. Chung and Graham \cite{chung-graham} later defined quasi-randomness for subsets of $\mathbb{Z}/N\mathbb{Z}$. Quasirandomness has been most thoroughly explored in the context of graphs, for which the reader should consult the excellent survey articles \cite{krivelevich-sudakov,komlos-simonovits}. The notions of uniformity in $\mathbb{Z}/N\mathbb{Z}$ and in $\mathbb{F}_p^n$ differ in little more than notation.\vspace{11pt}

\noindent As an example of uniformity/quasirandomness at work, and to get comfortable with the notation, let us prove that uniformity is more-or-less equivalent to a combinatorial condition involving $M(A)$, the number additive quadruples in $A$ (solutions to $a_1 + a_2 = a_3 + a_4$, $a_i \in A$). 

\begin{lemma}\label{add-quad-lem} Let $A \subseteq \mathbb{F}_p^n$ have cardinality $\alpha N$. 
\begin{enumerate}
\item Suppose that $A$ is $\eta$-uniform. Then $M(A) \leq (\alpha^4 + \eta^2 \alpha)N^3$.
\item Suppose that $M(A) \leq (\alpha^4 + \epsilon) N^3$. Then $A$ is $\epsilon^{1/4}$-uniform.
\end{enumerate}
\end{lemma}
\noindent\textit{Remark.} An easy application of the Cauchy-Schwarz inequality confirms that $M(A) \geq \alpha^4 N^3$, so this lemma concerns sets with close to the minimum number of additive quadruples.\vspace{11pt}

\proof The proof of this Lemma rests on the identity
\[ M(A) = N^{-1}\sum_{\xi} |\widehat{A}(\xi)|^4,\] which may be proved by observing that $M(A) = \sum_x (A \ast A)(x)^2$ and using Lemma \ref{basic-fourier-properties} (2) and (4). To prove (1), assume that $A$ is $\eta$-uniform, so that $|\widehat{A}(\xi)| \leq \eta N$ for all $\xi \neq 0$. Then we have
\begin{eqnarray*}
N M(A) & = & |\widehat{A}(0)|^4 + \sum_{\xi \neq 0} |\widehat{A}(\xi)|^4 \;\leq \; |A|^4 + \sup_{\xi \neq 0} |\widehat{A}(\xi)|^2\sum_{\xi} |\widehat{A}(\xi)|^2 \\ & = & \alpha^4 N^4 + \sup_{\xi \neq 0} |\widehat{A}(\xi)|^2 \cdot \alpha N^2 \; \leq \; (\alpha^4 + \eta^2 \alpha) N^4,\end{eqnarray*} as required.
To prove (2), assume that $M(A) \leq (\alpha^4 + \epsilon) N^3$. Then for any $\xi \neq 0$ one has
\[
|\widehat{A}(\xi)|^4 \; \leq \; \sum_{\xi} |\widehat{A}(\xi)|^4 - |\widehat{A}(0)|^4 \; = \; NM(A) - |A|^4 \; \leq \; \epsilon N^3,\]
which is what we wanted to prove.\endproof\vspace{11pt}

\noindent We observe that if $A \subseteq \mathbb{F}_p^n$, and if $H + g$ is a coset of some subspace $H \leq \mathbb{F}_p^n$, then there is a natural notion of what it means for $A$ to be $\eta$-regular relative to $H + g$.
Indeed we may define a set $\Agh \subseteq H$ by setting
\[ \Agh(x) \; = \; A(x + g)\] for $x \in H$. Since $H$ is a subgroup, it is isomorphic to $\mathbb{F}_p^m$ for some $m \leq n$ and it makes sense to talk about the Fourier transform on $H$. We say that $A$ is $\eta$-uniform on $H + g$ if $\Agh$ is $\eta$-uniform, considered as a subset of $H$.\vspace{11pt}

\noindent The key reason for uniformity being so important to us in the present survey is that it allows us to count solutions to certain linear equations in sets which are sufficiently uniform. Lemma \ref{add-quad-lem} was of course a rather special example of this (the linear equation being $a_1 + a_2 = a_3 + a_4$).  The next proposition illustrates this further.

\begin{proposition}\label{prop3.3}
Let $p$ be a prime and suppose that $A \subseteq \mathbb{F}_p^n$. Suppose that $\lambda_1,\dots,\lambda_k$, $k \geq 3$, are non-zero integers coprime to $p$. Let $H \leq \mathbb{F}_p^n$ be a subspace, and let $g_1,\dots,g_k \in \mathbb{F}_p^n$ satisfy $\sum \lambda_i g_i = 0$. Suppose that the density of $A$ on $H + g_i$ is $\alpha_i$, and that $A$ is $\eta$-uniform. Then $M$, the number of solutions to $\sum \lambda_i a_i = 0$ with $a_i \in H + g_i$ for $i = 1,\dots,k$, satisfies
\[ |M - \alpha_1 \dots \alpha_k |H|^{k-1}| \leq \eta^{k-1} (\alpha_1 \dots \alpha_k)^{1/k} |H|^{k-1}.\]
\end{proposition}
\proof With the notation introduced above we can write
\[ M = \sum_{\substack {h \in H \\ \sum \lambda_i h_i = 0 }} A_H^{+ x_1}(h_1) \dots A_H^{+x_k}(h_k).\]
This can be written in terms of the Fourier transform on $H$ as
\[ M = |H|^{-1}\sum_{\xi} \widehat{A_H^{+x_1}}(\lambda_1 \xi)\dots\widehat{A_H^{+x_k}}(\lambda_k \xi).\]
Separating off the term $\xi = 0$ and bounding the other term using H\"older's inequality, we get
\begin{eqnarray*}
|M - \alpha_1 \dots \alpha_k |H|^{k-1}| & \leq & |H|^{-1}\sum_{\xi \neq 0} |\widehat{A_H^{+x_1}}(\lambda_1 \xi)\dots\widehat{A_H^{+x_k}}(\lambda_k \xi)| \\ & \leq & \sup_{\xi \neq 0} |\prod_j \widehat{A_H^{+x_j}}(\lambda_j \xi)|^{1 - 2/k} \cdot \prod_{j} (\sum_{\xi} |\widehat{A_H^{+x_j}}(\lambda_j \xi)|^2)^{1/k}\\ & \leq & \eta^{k-2} (\alpha_1 \dots \alpha_k)^{1/k} |H|^{k-1}.
\end{eqnarray*}
This concludes the proof.\endproof\vspace{11pt}

\noindent Of particular importance to us will be two cases of the above with $k = 3$: $(\lambda_1,\lambda_2,\lambda_3) = (1,1,-2)$, which corresponds to arithmetic progressions of length 3, and $(\lambda_1,\lambda_2,\lambda_3) = (1,1,-1)$, corresponding to what are known as \textit{Schur triples} (solutions to $x + y = z$).\vspace{11pt}

\noindent A particularly nice feature of finite fields is that the notion of a set $A \subseteq \mathbb{F}_p^n$ being uniform is closely related to that set being well-distributed in cosets of codimension one hyperplanes. We will use this principle several times in the sequel, so let us state and prove a quantitative version of it now.

\begin{lemma}\label{uniformity-hyperplanes}
Suppose that $A \subseteq \mathbb{F}_p^n$ is a set of size $\alpha N$ \textup{(}$N = p^n$\textup{)} and that $A$ is not $\eta$-uniform, so that there is $\xi \neq 0$ with $|\widehat{A}(\xi)| > \eta N$. Let $H = \langle \xi \rangle^{\perp}$, and write $h(x) = H(x)/|H|$. Then
\begin{enumerate}
\item $\mathbb{E}(A \ast h(x)^2) \geq \alpha^2 + \eta^2$;
\item $\sup_x A \ast h(x) \geq \alpha + \frac{\eta^2}{\alpha}$;
\item $\sup_x A \ast h(x) \geq \alpha + \eta/2 $.
\end{enumerate}
\end{lemma}
\noindent\textit{Remark.} $A \ast h(x)$ is the density of $A$ on the coset $H + x$.\\
\proof To prove (1), observe that 
\begin{eqnarray*} N\sum_x A \ast h(x)^2 & = &  \sum_{\gamma} |\widehat{A}(\gamma)|^2 |\widehat{h}(\gamma)|^2 \\ & \geq & |\widehat{A}(0)|^2 |\widehat{h}(0)|^2 + |\widehat{A}(\xi)|^2 |\widehat{h}(\xi)|^2 \\& \geq & (\alpha^2 + \eta^2)N^2.\end{eqnarray*}
Statement (2) is a simple corollary of this:
\[ \alpha N \sup_x A \ast h(x) = \left(\sum_x A \ast h(x)\right) \sup_x A\ast h(x) \geq \sum_x A\ast h(x)^2.\]
Statement (3) is proved by working directly with the definition of $\widehat{A}(\xi)$. It leads to somewhat better qualitative bounds than (2).\vspace{11pt}

\noindent Let, then,  $H + x_j$, $j = 0,1,\dots,p-1$, be a complete set of cosets of $H$. Then \[
\widehat{A}(\xi) \; = \; \sum_j|A \cap H_j|\omega^j  \; = \; \sum_j a_j\omega^j,\] where $a_j = |A \cap H_j| - \alpha |H|$. Thus $\sum_j |a_j| \geq \eta N$. Observe, however, that $\sum a_j = 0$; it follows that $\sum_j |a_j| + a_j \geq \eta N$, and whence from the pigeonhole principle that $|a_j| + a_j \geq \eta N/p$ for some $j$. For such a $j$, we have $a_j \geq \eta N/2p$.\endproof\vspace{11pt}

\section{Roth's Theorem and the iteration method}\label{sec4}

\noindent Let us begin by recalling Problem \ref{prob1}.\vspace{11pt}

\noindent\textbf{Problem \ref{prob1}.} What is the cardinality of the largest subset of $\{1,\dots,N\}$ containing no three distinct elements $x,x+d,x+2d$ in arithmetic progression?\vspace{11pt}

\noindent This question was first raised by Erd\H{o}s and Tur\'an in 1936 \cite{erdos-turan}, and was addressed by
Klaus Roth \cite{roth}. Define $r_3(N)$ to be the answer to Problem \ref{prob1}. Roth proved that $r_3(N) \ll N/\log \log N$, a bound which was improved to $N(\log N)^{-c}$ independently by Heath-Brown \cite{heath-brown} and Szemer\'edi \cite{szemeredi-aps}, and then to $r_3(N) \ll N(\log\log N/\log N)^{1/2}$ by Bourgain \cite{bourgain2}. We are still a long way from a complete understanding of $r_3(N)$; the best known lower bound is Behrend's \cite{behrend} 1946 example showing that $r_3(N) \gg Ne^{-c\sqrt{\log N}}$.\vspace{11pt}

\noindent It is natural to define $r_3(G)$ for any group $G$ with no $2$-torsion (though see \cite{lev}). A particularly appealing case, which fits with the discussion of this article, is $G = \mathbb{F}_3^n$. In this case it turns out that the four proofs \cite{bourgain2,heath-brown,roth,szemeredi-aps} can all be adapted to give the following result.

\begin{theorem}\label{finite-field-roth}
We have $r_3(\mathbb{F}_3^n) \ll N/\log N$ \textup{(}= $O(3^n/n)$\textup{)}.
\end{theorem}
\noindent In fact, all four proofs look the same in the finite field setting. Roth's proof was adapted to the finite field setting by Meshulam \cite{meshulam} and the argument we give to prove Theorem \ref{finite-field-roth} is the same as his.\vspace{11pt}

\noindent There are two key ingredients. The first is a special case of Proposition \ref{prop3.3}, asserting that if $A$ is sufficiently uniform then we can count solutions to the equation $a_1 + a_2 = 2a_3$ (that is, arithmetic progressions of length three).

\begin{lemma}\label{ap3-lemma}
Suppose that $A \subseteq \mathbb{F}_3^n$ has cardinality $\alpha N$, and that $A$ is $\eta$-uniform. Then there are at least $(\alpha^3 - \eta \alpha)N^2$ solutions to the equation $a_1 + a_2 = 2a_3$ with $a_i \in A$. In particular if $\eta = \alpha^2/2$ and $N > 2/\alpha^2$ then $A$ contains a 3-term AP $(x,x+d,x+2d)$ with $d \neq 0$.
\end{lemma}
\proof The first part is just a matter of setting $H = \mathbb{F}_3^n$ and $(\lambda_1,\lambda_2,\lambda_3) = (1,1,-2)$ in Proposition \ref{prop3.3}.
To verify the second statement, one must simply check that if $\eta  = \alpha^2 /2$ and $N > 2 /\alpha^2$ then $(\alpha^3 - \eta \alpha)N^2$ is greater than $\alpha N$, the number of ``trivial'' 3-term APs $(x,x,x)$ in $A$.\endproof\vspace{11pt}

\noindent The second key ingredient is Lemma \ref{uniformity-hyperplanes} (3), which asserts that if $A$ is not $\eta$-uniform then it has increased density on some coset of a hyperplane.
In combination with Lemma \ref{ap3-lemma} this leads naturally to an \textit{iterative} method for proving Theorem \ref{finite-field-roth}. \vspace{11pt}

\noindent\textit{Proof of Theorem \ref{finite-field-roth}.} Set $A_0 = A$, $H_0 = \mathbb{F}_3^n$, $\alpha_0 = \alpha$. For each $i = 0,1,\dots$ we perform the following algorithm:
\begin{itemize}
\item If $A_i$ is $\alpha_i^2/2$-uniform then \texttt{STOP}. 
\item Otherwise by Lemma \ref{uniformity-hyperplanes} find a hyperplane $H_{i+1} \leq H_i$ and an $x \in H_i$ such that $|A_i \cap (H_{i+1} + x)| > (\alpha_i + \alpha_i^2/4)|H_{i+1}|$. Now set $A_{i+1} = (A_i - x) \cap H_{i+1}$ and set $\alpha_{i+1} = |A_i \cap (H_{i+1} + x)|/|H_{i+1}|$.
\end{itemize}
Note that if $A_i$ contains a 3-term AP then so does $A$.\vspace{11pt}

\noindent The algorithm cannot be repeated forever, since the sequence $(\alpha_i)_{i=1}^{\infty}$ satisfies $\alpha_0 = \alpha$ and $\alpha_{i+1} \geq \alpha_i + \alpha_i^2/4$ then $\alpha_i > 1$ for $i > 50/\alpha$. Thus we reach a \texttt{STOP} at step $K$ of the algorithm, for some $K < 50/\alpha$. At this stage, $A_K$ is $\alpha_K^2/2$-uniform. If in addition $|H_K| > 2/\alpha_K^2$ then, by Lemma \ref{ap3-lemma}, $A_K$ contains a 3-term AP. Since $|H_K| = 3^{-K}N > 3^{-50/\alpha}N$ and $\alpha_K > \alpha$, we see that the original set $A$ contains a 3-term AP if $\alpha > C/\log N$ for some $C$.\endproof\vspace{11pt}

\noindent We call the above an \textit{iteration argument} for obvious reasons. We will encounter several such arguments in this survey, so let us take the opportunity to look at the important features of it.\vspace{11pt}

\noindent Our concern was with certain \textit{configurations} $\mathbf{Config}$, which in this section were the three-term arithmetic progressions $(x,x+d,x+2d)$, $d \neq 0$.\vspace{11pt}

\noindent A key feature of the argument was a collection $\mathbf{Struct}$ of \textit{structures}, which in this case was the collection of all cosets of subspaces of $\mathbb{F}_3^n$. There was also some measure of the \textit{complexity} $\omega(S)$ of a given structure $S \in \mathbf{Struct}$, this being the codimension of the subspace. For a given set $A \subseteq \mathbb{F}_3^n$  and for any $S \in \mathbf{Struct}$ there was a notion of the \textit{density} $\delta_S(A)$ of $A$ relative to $S$. Finally, there was a \textit{norm} $\Vert \cdot \Vert_S$ on functions $f : S \rightarrow [-1,1]$, for any $S \in \mathbf{Struct}$ (in the example above, this was the $L^{\infty}$ norm of the Fourier transform of $f$, regarded as a function on $\dot{S}$). This we used to define a notion of \textit{uniformity} relative to some $S \in \mathbf{Struct}$; a set was $\eta$-uniform if $\Vert A - \delta_S(A) \Vert_S \leq \eta$. \vspace{11pt}

\noindent The ``iteration step'' of Roth's argument can be presented in the following way.\vspace{11pt}

\noindent Let $S \in \mathbf{Struct}$, and let $A \subseteq \mathbb{F}_3^n$ be a set with $\delta_S(A) = \alpha$. Then one of the following three alteratives holds:
\begin{enumerate}
\item (generalised von Neumann theorem\footnote{This term is one that Tao and I are trying to popularize to emphasise the connection with results in ergodic theory such as \cite[Lemma 3.1]{FKO}. Such results tend to be established using several applications of the Cauchy-Schwarz inequality -- see for example \cite[\S 5]{green-tao-primes}. The phrase ``key lemma'' was used for a related concept in the theory of graph regularity in the excellent survey \cite{komlos-simonovits}: now the more descriptive term ``counting lemma'' is popular (cf. \cite{gowers-hyper,green-reg,nagle-rodl-schacht}).}) $\Vert A - \delta_S(A) \Vert_S \leq \alpha^2/2$, in which case $A$ contains some $c \in \mathbf{Config}$;
\item (density increment) $\Vert A - \delta_S(A) \Vert_S > \alpha^2/2$, in which case we may find $S' \in \mathbf{Struct}$, $\omega(S') \leq \omega(S) + 1$, such that $\delta_{S'}(A) \geq \delta_S(A) + \alpha^2/4$;
\item (endpoint) $|S| < 2/\alpha^2$.
\end{enumerate}

\noindent Several subsequent arguments will have the same general form, with different notions of \textbf{Struct}, \textbf{Config}, $\omega$ and $\Vert \cdot \Vert_S$. The choice of \textbf{Struct} and, perhaps more importantly, of the norm $\Vert \cdot \Vert_S$ is vitally important. $\Vert \cdot \Vert_S$ must be ``strong'' enough for us to be able to prove a von Neumann theorem, yet ``weak'' enough that one may obtain a density increment.\vspace{11pt}

\noindent To conclude this section, let use return to the question of estimating $r_3(\mathbb{F}_3^n)$, which I regard as a very interesting one. It seems to dramatically expose our lack of understanding of 3-term arithmetic progressions. There does not seem to be an analogue of Behrend's example in the finite field setting (Behrend's construction makes important use of convexity in $\mathbb{R}^n$). The best known lower bounds on $r_3(\mathbb{F}_3^n)$ come from design theory, where a set in $\mathbb{F}_3^n$ with no 3-term AP is known as a \textit{cap}. Write $f(n)$ for the cardinality of the largest cap in $\mathbb{F}_3^n$. In \cite{edel} one finds the estimate
\[ \mu(3) := \limsup_{n \rightarrow \infty} \frac{\log_3(f(n))}{n} \geq 0.724851,\] which seems to be the best known.
In that paper it is stated as an interesting research problem to determine if $\mu(3) = 1$. I believe that this is not so.

\begin{conjecture}
$\mu(3) < 1$. That is, there is an absolute constant $\delta > 0$ such that $r_3(\mathbb{F}_3^n) \leq (3 - \delta)^n$.
\end{conjecture}

\noindent I would expect any methods used to make progress on this conjecture to assist with the Problem \ref{prob1}. At present the best known bound is that given in Theorem \ref{finite-field-roth}.

\section{Right-angled triangles - an argument of Shkredov}\label{sec5}
\noindent In this section we write $V_n = \mathbb{F}_2^n$, and $N = |V_n| = 2^n$.

\noindent We are concerned with a sort of two-dimensional generalisation of Problem \ref{prob1}:

\begin{problem} What is $r_{\angle}(N)$, the cardinality of the largest subset of $\{1,\dots,N\} \times \{1,\dots,N\}$ containing no \emph{corner} $((x,y), (x + d,y), (x, y+d))$, $d \neq 0$?\end{problem}

\noindent Ajtai and Szemer\'edi \cite{ajtai-szemeredi} proved that $r_{\angle}(N) = o(N)$, and various subsequent authors \cite{solymosi,vu} have obtained explicit bounds of the shape $r_{\angle}(N) \ll N/(\log_* N)^c$. Here $\log_* N$ is the number of times one must take the logarithm of $N$ in order to produce a number less than 2. \vspace{11pt}

\noindent Very recently Shkredov \cite{shkredov} produced the first ``sensible'' bound \[ r_{\angle}(N) \ll N/(\log \log \log N)^c.\]  In this section we give the finite field version of his argument, in which the details are greatly simplified.\vspace{11pt}

\noindent Let $G$ be an abelian group of size $N$, and consider the collection of corners in $G \times G$, by which we mean triples $((x,y), (x+d,y), (x,y+d))$, $d \neq 0$. Write $r_{\angle}(G)$ for the cardinality of the largest set $A \subseteq G \times G$ which does not contain any corner.

\begin{theorem}[Shkredov] \label{shkredov-thm}We have $r_{\angle}(\mathbb{F}_2^n) \ll N^2/(\log \log N)^{1/25}$.
\end{theorem}

\noindent It is natural to try and use the iteration method, in the form outlined in the previous section. The most na\"{\i}ve attempt at doing this would involve taking $\mathbf{Struct}$ to be the set of cosets of products $H \times H$, where $H \leq \mathbb{F}_2^n$ is a subspace, and the definition of uniformity to be much the same as before. The notion of having no large Fourier coefficients makes perfect sense in $H \times H$. Unfortunately, however, this notion of uniformity is not subtle enough to give good control on the number of corners, essentially because it does not ``see'' the coordinate structure of $H \times H$. The following example is instructive:

\begin{example}\label{ex1}
Let $B$ be a random (and hence highly uniform) subset of $V_n$ with cardinality $\beta N$, and let $A \subseteq \mathbb{F}_2^n \times \mathbb{F}_2^n$ be the set $B \times B$. Then $A$ is also highly uniform. The density of $A$ is $\alpha = \beta^2$. A corner in $A$ corresponds to a quadruple of points $(x,x+d,y,y+d) \in B^4$, and we know from Lemma \ref{add-quad-lem} that there are roughly $\beta^4 N^3 = \alpha^2 N^3$ such configurations. If $A$ were truly random, however, it would have more like $\alpha^3 N^3$ corners.
\end{example}

\noindent The next idea, then, might be to 
define a somewhat finer notion of uniformity which respects the coordinate structure somewhat more. Using Proposition \ref{add-quad-lem} as a guide, we might define $A \subseteq \mathbb{F}_2^n \times \mathbb{F}_2^n$ to be \textit{rectilinearly $\eta$-uniform} if the number of configurations $((x,y),(x+d,y),(x,y+e),(x+e,y+e))$ in $A^4$ is at most $(\alpha^4 + \eta) N^4$. Such a notion does, as we will see, give some control on the number of corners in $A$. Unfortunately passing to a new structure $S \in \mathbf{Struct}$ on which the density increases is now problematic.\vspace{11pt}

\noindent To see why, consider again example \eqref{ex1}. It is easy to see that $A$ fails to be rectilinearly uniform, but there is no product set $(H + x) \times (H' + x')$, $H, H'$ large subspaces of $\mathbb{F}_2^n$, on which the density of $A$ increases markedly. \vspace{11pt}

\noindent Note, however, that in this example there \textit{is} at least some structure on which the density of $A$ increases, and that is the product set $B \times B$ (of course, the density of $A$ on this set is one). This behaviour is more-or-less typical: if a set $A \subseteq \mathbb{F}_2^n \times \mathbb{F}_2^n$ has substantially more than $\alpha^4 N^4$ rectangles then it has increased density on some product $B_1 \times B_2$. This can be proved graph-theoretically by associating to $A$ the bipartite graph $\Gamma_A$ with vertex sets two copies of $\mathbb{F}_2^n$, an edge $xy$ being deemed to lie in $\Gamma_A$ precisely if $(x,y) \in A$. A rectangle in $A$ then corresponds to a copy of $C_4$ in $\Gamma_A$, and we are reduced to showing that if $\Gamma_A$ has substantially more than $\alpha^4 N^4$ copies of $C_4$ then there are large vertex sets $B_1,B_2$ such that the edge density of $\Gamma_A$ restricted to $B_1 \cup B_2$ is much greater than $\alpha$. Shkredov in effect provides a spectral proof of this statement, which is in the spirit of \cite{frieze-kannan}. A purely combinatorial proof is more traditional, and somewhat simpler -- the details may be found in \cite{sup-doc-1}.\vspace{11pt}

\noindent The discussion of the last paragraph might suggest that we should enlarge $\mathbf{Struct}$ to include all translates of products $B \times B$. This turns out to be too much of a compromise -- one cannot establish a useful generalised von Neumann theorem. \vspace{11pt}

\noindent The above discussions motivate Shkredov's main advance, which is an appropriate definition of $\mathbf{Struct}$. The definition depends on the global density $\alpha$ of $A$, a feature which has no analogue in other applications of the iterative method discussed in this paper. 

\begin{definition}\label{struct-def}
Let $\alpha > 0$. Define $\mathbf{Struct}_{\alpha}$ to consist of all translates of product sets $S = E_1 \times E_2$, where $E_1,E_2$ are subsets of some $H \leq \mathbb{F}_2^n$, $|E_i| = \beta_i |H|$ and $E_i$ is a $(2^{-36}\beta_1^{12}\beta_2^{12} \alpha^{36})$-uniform subset of $H$ for $i = 1,2$. 
\end{definition}

\begin{definition}\label{rect-norm} Suppose that $S = E_1 \times E_2$ is a product set,  and that $f : S \rightarrow [-1,1]$ is a function. Then we define the \emph{rectangle norm} of $f$, $\Vert f \Vert_{S}$ by
\[ \Vert f \Vert_{S}^4 = \mathbb{E}\left( f(x,y)f(x',y)f(x,y')f(x',y') | x,x' \in E_1, y,y' \in E_2\right).\] 
\end{definition}

\noindent It is not totally obvious that $\Vert \cdot \Vert_S$ is a norm, but this is in fact the case.
Let us now look at how the argument fits together, starting with a generalised von Neumann theorem.

\begin{proposition}[Generalised von Neumann]\label{gvn}
Let $S \in \mathbf{Struct}_{\alpha}$, so that $\mathcal{S} = E_1 \times E_2$ be a product set, where $E_1,E_2 \subseteq H$, $|E_i| = \beta_i |H|$ and $E_i$ is $(2^{-36}\beta_1^{12}\beta_2^{12} \alpha^{36})$-uniform for $i= 1,2$. Let $A \subseteq S$ be a set with $\delta_{S}(A) \geq  \alpha$. Suppose that  and that $\Vert A - \delta_S(A)\Vert_{\Box}^4 \leq 2^{-8}\alpha^{12}$. Then $A$ has at least $\alpha^3 \beta_1^2\beta_2^2 N^3/2$ corners.
\end{proposition}
\noindent The proof of this statement involves a number of applications of Cauchy-Schwarz.\vspace{11pt}

\noindent To complement the generalised von Neumann theorem, we must establish a density increment result.
The following can be obtained by simple graph theory (or alternatively by spectral methods, as done in \cite{shkredov}).

\begin{proposition}[Density increment on a product set]\label{prop5.7} Let $S = E_1 \times E_2$ be a product set, and suppose that $A \subseteq S$ has $\delta_{S}(A) = \alpha$ and $\Vert A - \alpha \Vert_{S}^4 \geq \eta$. Then there are sets $F_i \subseteq E_i$ with $|F_i| \geq 2^{-8}\eta|E_i|$ such that the density of $A$ on $S' = F_1 \times F_2$ satisfies $\delta_{S'}(A) \geq \alpha + 2^{-14}\eta^2$.
\end{proposition}
\noindent\textit{Remark.} There is no need to assume that the sets $E_1,E_2$ are uniform in this proposition.

\noindent Proposition \ref{prop5.7} has a significant deficiency, which means that it cannot be used in combination with Proposition \ref{gvn} to provide an iterative proof of Theorem \ref{shkredov-thm}. This is that the sets $F_1,F_2$ which it outputs need not be uniform, and so it is quite possible that $S' \notin \mathbf{Struct}_{\alpha}$. The following further result is required.

\begin{proposition}[Uniformising a product set]\label{prop5.8}
Let $\alpha,\tau,\sigma \in (0,1)$ be parameters, and let $S' = F_1 \times F_2$ be a product set in $W \times W$ with $|F_i| = \delta_i N$. Suppose that $A \subseteq S'$ is a set with $\delta_{S'}(A) = \alpha + \tau$, and that
\begin{equation}\label{M-cond}
|W| \geq \exp(16 \sigma^{-2}\delta^{-1}\tau^{-1}).
\end{equation}
Then there is a subspace $W' \subseteq W$, $\dim W' \geq \dim W - 8\sigma^{-2}\delta^{-1}\tau^{-1}$ and $t_1,t_2 \in W$ such that if $E'_1 = (F_1 - t_1) \cap W'$, $E'_2 = (F_2 - t_2) \cap W'$ and $S'' = E'_1 \times E'_2$ then
\begin{enumerate}
\item $|S''| \geq \delta_1\delta_2 \tau |W'|^2/2$;
\item $E'_1, E'_2$ are $2\sigma$-uniform as subsets of $W'$;
\item $\delta_{S''}(A - (t_1,t_2)) \geq \alpha + \tau/8$.
\end{enumerate}
\end{proposition}
\noindent The proof of this theorem also proceeds by a version of the iterative method, and in this sense Skhredov's argument is a sort of double iteration method. The most important content of the proposition is that if $S \subseteq W \times W$ then we may pass to a translate of $W' \times W'$ on which $S$ looks uniform, where $W' \leq W$ is a subspace of somewhat large codimension. If this really was our only aim, then we could proceed as follows. Either $S$ is already uniform, or else $S$ has a large Fourier coefficient $\xi$. In the latter case, $S$ has increased density on some translate of $\xi^{\perp}$, by Lemma \ref{uniformity-hyperplanes} (2). $\xi^{\perp}$ obviously contains a set of the form $W' \times W'$, with $W'$ having codimension at most two. Now simply iterate the argument. \vspace{11pt}

\noindent The one further issue is that we also need to keep control of the density of $A$, which sits inside $S$. To achieve this it is necessary to \textit{partition} $W \times W$ into pieces which are translates of products $W' \times W'$, such that $S$ is uniform on almost all of them. By a simple pigeonhole argument there must be some piece on which $S$ is uniform, and on which the relative density of $A$ is still quite large. Note that the subspaces $W'$ need not be the same for each piece; this is important from the point of view of obtaining bounds, or else one runs into examples such as that in \S 9 of \cite{green-reg}.\vspace{11pt}

\noindent To get this decomposition into pieces one uses the iterative argument with one small modification. At the $j$th stage of the iteration we will have a collection $\mathcal{C_j}$ of pieces, each being a translate of some product $W' \times W'$. If $c \in \mathcal{C}_j$, write $\delta(c)$ for the relative density of $S$ on the piece $c$. Our previous proposal was to ensure that $\sup_{c \in \mathcal{C}_j} \delta(c)$ increases at each step of the iteration, this idea having served us well in the past.
What one does instead is to increase the $L^2$ average $\mathbb{E}(\delta(c)^2 | c \in \mathcal{C}_j)$. This can be accomplised by using Lemma \ref{uniformity-hyperplanes} (1).\vspace{11pt}

\noindent Propositions \ref{prop5.7} and \ref{prop5.8} together give the requisite density increment result to go with the generalised von Neumann theorem of Proposition \ref{gvn}. Thus we can employ an iteration argument. Working out the bounds gives Theorem \ref{shkredov-thm}.

\section{Progressions of Length Four}\label{sec6}

\noindent In this section we give another example of the iterative method at work.

\begin{problem} Estimate $r_4(N)$, the cardinality of the largest subset of $\{1,\dots,N\}$ containing no four distinct elements $x,x+d,x+2d,x+3d$ in arithmetic progression?
\end{problem}
\noindent This question was, like Problem \ref{prob1}, raised by Erd\H{o}s and Tur\'an in 1936. Szemer\'edi \cite{szemeredi-4} was the first to show that $r_4(N) = o(N)$. It was not until as recently as 1998 that the first ``sensible'' upper bound, $r_4(N) \ll N(\log \log N)^{-c}$, was provided by Gowers \cite{gowers-szem4}. Gowers' argument was iterative, like the arguments of \S \ref{sec4} and \ref{sec5}. \vspace{11pt}

\noindent Of course, one can define $r_4(G)$ for any abelian group $G$ of size $N$. Recently, T.Tao and the author \cite{green-tao-ap4s} studied the case $G = \mathbb{F}_5^n$, starting from Gowers' work. Certain features of \cite{gowers-szem4} become rather simpler in this setting, and we were able to run the iterative method quite efficiently, obtaining the following theorem.

\begin{theorem}[G.--Tao] $r_4(\mathbb{F}_5^n) \ll N(\log N)^{-c}$ for some $c > 0$.
\end{theorem}

\noindent Write $\mathbf{Config}$ for the collection of all four-term progressions in $\mathbb{F}_5^n$. Any hope of proving a generalised von Neumann theorem with the same uniformity norm that we used in \S \ref{sec4} is dashed by the following example:
\begin{example}[Gowers; Furstenberg-Weiss]
There is a set $A \subseteq \mathbb{F}_5^n$ with density $1/5$, which is highly uniform, but which does not contain roughly $5^{-4} N^2$ four-term arithmetic progressions.
\end{example}
\proof Let $A = \{x \in \mathbb{F}_5^n : x^Tx = 0\}$. Then $A$ certainly has density approximately $1/5$. To see that $A$ is highly uniform, write
\[ \widehat{A}(\xi) = \frac{1}{5}\sum_{\lambda \in \mathbb{F}_5} \sum_{x \in \mathbb{F}_5^n} \omega^{\lambda x^T x - \xi^T x} = \frac{1}{5}\sum_{\lambda} \prod_{j=1}^n \omega^{\lambda x_j^2 - \xi_j x_j}.\]
If $\lambda \neq 0$ then each term in the product has magnitude $\sqrt{5}$, giving a total contribution of $5^{n/2}$; if $\lambda = 0$ then, provided $\xi \neq 0$, at least one term in the product vanishes. It follows that $\sup_{\xi \neq 0} |\widehat{A}(\xi)| \leq 5^{n/2} = \sqrt{N}$.\vspace{11pt}

\noindent However, $A$ has roughly $5^{-3}N^2$ progressions of length four. Indeed, since $A$ is so highly uniform we know from Proposition \ref{prop3.3} that it contains roughly this many progressions of length \textit{three}. However if $x$, $x+d$ and $x+2d$ all lie in $A$ then $x + 3d \in A$ automatically, in view of the easily verified identity
\begin{equation}\boxeq x^Tx - 3(x+d)^T(x+d) + 3(x+2d)^T(x + 2d) - (x+3d)^T(x+3d) = 0.\end{equation}
\noindent\textit{Remark.} Gowers has shown us an example of a subset of $\mathbb{Z}/N\mathbb{Z}$ which is uniform and has density $\alpha$, but has many \textit{fewer} than $\alpha^4 N^2$ four-term arithmetic progressions. \vspace{11pt}

\noindent Similar examples can be constructed using any quadratic form $q(x) = x^T M x + r^Tx + b$ in place of $x^T x$. Remarkably, there are essentially no other examples. We shall formalise this statement in what follows.

\begin{definition}[Gowers norm]
Let $f : \mathbb{F}_5^n \rightarrow [-1,1]$ be a function. Then the Gowers $U^3$-norm of $f$, $\Vert f \Vert_{U^3}$, is defined by $\Vert f \Vert_{U^3}^8 = $
\begin{equation}\label{gowu3}  \mathbb{E}(f(x)f(x+a)f(x+b)f(x+c)f(x+a+b)f(x+a+c)f(x+b+c)f(x+a+b+c)|x,a,b,c).\end{equation}
\end{definition}
\noindent Again, it is not completely obvious that $\Vert \cdot \Vert_{U^3}$ is a norm, but this is not to hard to show. The following result is due to Gowers \cite{gowers-szem4}. As with the other generalised von Neumann theorems we have mentioned, the proof involves several applications of the Cauchy-Schwarz inequality.

\begin{theorem}(Generalised Von Neumann theorem)\label{gow-thm}
Suppose that $A \subseteq \mathbb{F}_5^n$ has density $\alpha$, and that $\Vert A - \alpha \Vert_{U^3} \leq \delta^4/100$. Then $A$ has at least $\alpha^4 N^2/2$ four-term arithmetic progressions.
\end{theorem}

\noindent The next theorem is proved in \cite{green-tao-ap4s} by adding a single new idea, the so-called ``symmetry argument'', to the ideas of Gowers \cite{gowers-szem4}. This theorem clarifies the sense in which the ``quadratic'' examples of Furstenberg and Weiss are in a sense the only ones:

\begin{theorem}[Gowers; G.--Tao]\label{gow-thm2} Suppose that $\Vert A - \alpha \Vert_{U^3} \geq \delta$. Then $A$ has quadratic bias, meaning that there is some quadratic form $q(x) = x^T M x + r^T x + b$ such that $A$ has density at least $\alpha +  C(\delta)$ on the zero set $S = \{x : q(x) = 0\}$.
\end{theorem}

\noindent The reader may note that these two theorems do not, in their present incarnations, dovetail together to give an iteration argument because there is no natural definition of $\mathbf{Struct}$. Roughly speaking, Gowers took $\mathbf{Struct}$ to be the collection of translates of subspaces of $\mathbb{F}_5^n$. It is possible to deduce from the conclusion of Theorem \ref{gow-thm2} that $A$ has increased density on some $S \in \mathbf{Struct}$, but unfortunately the codimension $\omega(S)$ might be exceedingly large (perhaps $n - n^{1/100}$). This does not, then, lead to a very efficient iterative argument.\vspace{11pt}

\noindent In \cite{green-tao-ap4s} a much less appetising approach is forced to work, which leads to superior bounds. Roughly, this involves taking $\mathbf{Struct}$ to be the collection of all \textit{quadratic submanifolds}, our name for an intersection 
\[ S = \bigcap_{j=1}^k \{x : q_j(x) = 0\},\]
where $q_j(x) = x^T M_j x + r_j^T x + b_j$ are quadratic forms. The ``roughly'' is quite important. We must in fact assume that $S$ is ``generic'', meaning that the matrices $M_j$ are not too linearly dependent. In practise this means that they satisfy a \textit{rank condition} such as $\mbox{rk}(\lambda_1 M_1 + \dots + \lambda_k M_k) \geq 10k$ for all possible choices of scalars $\lambda_j \in \mathbb{F}_5$. We also allow our quadratic forms to be defined only on a subspace $W \leq \mathbb{F}_5^n$ of not-too-large codimension. This is because of the very useful observation that an arbitrary quadratic submanifold can be made generic after passing to an appropriate subspace $W$.\vspace{11pt}

\noindent Generalising the notion of Gowers $U^3$-norm to such a setting is straightforward; in fact the definition is the same except that the expectation in \eqref{gowu3} is taken over $S$. Proving an analogue of Theorem \ref{gow-thm} is substantially more involved, but it is possible and reads as follows.

\begin{theorem}\label{gow-thm22} Let $S \in \mathbf{Struct}$. That is to say, $S$ is a generic quadratic submanifold in some $W \leq \mathbb{F}_5^n$, this being the zero set of some $k$ quadratic forms $q_1,\dots,q_k$ on $W$. Then $S$ has approximately $|S|^3/|W|$ four-term arithmetic progressions.
\begin{enumerate}
\item \textup{(Generalised von Neumann theorem)} Suppose that $A \subseteq S$ has density $\alpha$, and that $\Vert A - \alpha \Vert_{U^3(S)} \leq \alpha^{20}$. Then $A$ has at least $\alpha^4 |S|^3/2|W|$ four-term arithmetic progressions.
\item \textup{(Gowers-type inverse theorem)} Suppose that $\Vert A - \alpha \Vert_{U^3(S)} \geq \delta$. Then $A$ has quadratic bias, meaning that there is some quadratic form $q_{k+1}(x) = x^T M_{k+1} x + r_{k+1}^T x + b_{k+1}$ such that $A$ has density at least $\alpha +  C(\delta)$ on the set $S \cap \{x : q_{k+1}(x) = 0\}$.
\end{enumerate}
\end{theorem}

\noindent A key feature of the theorem is that the density increment $C(\delta)$ is independent of the number of quadratic forms $k$. The proof of the theorem is long and somewhat difficult and occupies the bulk of \cite{green-tao-ap4s}.\vspace{11pt}

\noindent Theorem \ref{gow-thm22} of course allows one to set up an iteration scheme. If $A \subseteq \mathbb{F}_5^n$ is a set with density $\alpha$ which contains no four-term progressions, then one may find a sequence \[ \mathbb{F}_5^n = S_0 \supseteq S_1 \supseteq S_2 \supseteq \dots \] of generic quadratic manifolds, defined on subspaces \[ \mathbb{F}_5^n = W_0 \geq W_1 \geq W_2 \geq \dots \] such that the density of $A$ on $S_j$ is at least $\alpha + jC(\alpha^{20})$. This leads to a contradiction after $C(\alpha^{20})^{-1}$ iterations.\vspace{11pt}

\noindent Unfortunately, this still leads to a bound of the shape $r_4(\mathbb{F}_5^n) \ll N(\log \log N)^{-c}$, since we have only been able to establish Theorem \ref{gow-thm22} with a function $C(\delta)$ which behaves like $\exp(-\delta^{-B})$, and this results in a very large number of iterations. We conjecture that a better bound holds, but we cannot prove this even in the less general context of Theorem \ref{gow-thm2}. I regard this as one of the key open questions in this area of arithmetic combinatorics.

\begin{conjecture}[Polynomial Gowers Inverse Conjecture]\label{quad-inverse}
Let $f : \mathbb{F}_5^n \rightarrow [-1,1]$ be a function with $\mathbb{E}f = 0$. Suppose that $\Vert f \Vert_{U^3} \geq \delta$. Then there is a quadratic form $q$ on $\mathbb{F}_5^n$ such that
\[ |\mathbb{E} f(x) \omega^{q(x)}| \gg \delta^{C},\] for some absolute constant $C$.
\end{conjecture}

\noindent We do know this with $\delta^{C}$ replaced by a function of exponential type. An affirmative answer to the PGI conjecture would be implied by an affirmative answer to the Polynomial Freiman-Ruzsa conjecture (PFR), which is discussed in some detail in \S \ref{sec10}.\vspace{11pt}

\noindent Fortunately, for the purposes of obtaining a bound on $r_4(\mathbb{F}_5^n)$ one can get by with a weaker conclusion in Theorems \ref{gow-thm} and \ref{gow-thm2}. In Theorem \ref{gow-thm}, one can obtain a ``polynomial'' density increment, leading to a much shorter iterative process, by passing to a set of the form $\{x : q(x) = 0\} \cap (W + t)$, where $W \subseteq \mathbb{F}_5^n$ is a subspace. One can allow the codimension of $W$ to be a power of $\alpha^{-1}$, which is just as well since this is the best bound we have.

\section{Szemer\'edi Regularity in Groups}\label{sec7}

\noindent The object of this section is to state some results and open problems from  \cite{green-reg}. The results are slightly different from those in the previous section in the problem addressed is not quite of ``Szemer\'edi type''. However what we discuss here is certainly in a similar spirit, being concerned with solutions of linear equations in sets of integers, and can furthermore be interpreted as an application of the iteration method.\vspace{11pt}

\noindent We will be somewhat brief: more details can of course be found in the paper \cite{green-reg} itself, which is written from a viewpoint rather similar to that of the present survey.\vspace{11pt}

\noindent Szemer\'edi's regularity lemma is a famous result in graph theory. It can be regarded as structure theorem for all graphs, in the sense that it shows that one can decompose a completely arbitrary graph into a bounded number of pieces, almost all of which are pseudorandom. There are many excellent articles on this topic -- see for example \cite{komlos-simonovits}.\vspace{11pt}

\noindent One consequence of Szemer\'edi's regularity lemma is the following interesting result.\footnote{We have not attributed this result, as it is not clear to us where it was first stated. A slightly weaker result was obtained by Ruzsa and Szemer\'edi in 1976 \cite{ruz-szem}. The result is also well-known in the literature concerning ``property testing'': see, for example, \cite{alon-property-testing}.}
\begin{theorem} \label{ruz-szem-thm} Let $\Gamma$ be a graph on $N$ vertices, and suppose that one must remove $\delta N^2$ edges from $G$ in order to destroy all triangles in $\Gamma$. Then $\Gamma$ has at least $C_1(\delta) N^3$ triangles, for some $C(\delta) > 0$.
\end{theorem}
\noindent Put another way, if a graph is almost triangle-free (i.e. contains few triangles) then it can be made truly triangle-free by the removal of a small number of edges.\vspace{11pt}

\noindent Our investigations in \cite{green-reg} were motivated by an ``arithmetic'' question related to the above theorem.
\begin{theorem}[See \cite{green-reg}] \label{almost-sum-free} Let $G$ be an abelian group of size $N$, and suppose that $A \subseteq G$ is a set. Suppose that one must remove $\delta N$ elements from $A$ in order to create a sum-free set \textup{(}that is, a set with no solutions to $x + y = z$\textup{)}. Then $A$ has at least $C_2(\delta) N^2$ Schur triples \textup{(}triples $(x,y,z)$ for which $x + y = z$\textup{)}.\end{theorem}
\noindent This result may be regarded as a structure theorem for sets which are almost sum-free; they can be made truly sum-free by the removal of a few elements.
\vspace{11pt}

\noindent This theorem is deduced from a result which we call a Szemer\'edi-type regularity lemma for abelian groups. This result is a perfect example for the present survey, since in the context of a general abelian group it requires substantial preparation to even state the result. When $G = \mathbb{F}_2^n$, however, things are much easier.

\begin{theorem}[Regularity lemma for $\mathbb{F}_2^n$]\label{reg-thm}
Let $A \subseteq \mathbb{F}_2^n$ be a set, and let $\epsilon > 0$ be a parameter. Then there is a subspace $H \subseteq \mathbb{F}_2^n$ with codimension at most $M(\epsilon)$, and such that $A$ is $\epsilon$-uniform on at least a proportion $1 - \epsilon$ of the cosets of $H$.
\end{theorem}

\noindent Let us say, for the rest of this section, that $A$ is $\epsilon$-regular relative to $H$ if it satisfies the conclusion of this theorem.\vspace{11pt}

\noindent Let us sketch the deduction of Theorem \ref{almost-sum-free} from Theorem \ref{reg-thm}. Suppose that $A \subseteq \mathbb{F}_2^n$ is a set with the property that one must remove at least $\delta N$ elements from $A$ to leave a set which is sum-free. Apply Theorem \ref{reg-thm} with $\epsilon = (\delta/10)^3$, giving a subspace $H$ of codimension at most $M(\epsilon)$ such that $A$ is $\epsilon$-uniform for a proportion at least $1 - \epsilon$ of the cosets of $H$. For each coset $H + x$, we ask two questions:
\begin{itemize}
\item Is $A$ $\epsilon$-uniform on $H + x$?
\item Is the density of $A$ on $H + x$ at least $(2\epsilon)^{1/3}$?
\end{itemize}
If the answer to either of these questions is \textit{no} then we simply remove all of $A \cap (H + x)$ from $A$. Let the set remaining after we have asked the above questions for all cosets $H + x$ be called $A'$. It is easy to see that \[ |A'| > |A| - 10 \epsilon^{1/3} N = |A| - \delta N.\]
We claim that $A'$ is sum-free. Indeed, were it not there would be $x_1,x_2,x_3$ with $x_1 + x_2 = x_3$, such that $A$ is $\epsilon$-uniform and has density $\alpha_i \geq (2\epsilon)^{1/3}$ on each $H + x_i$. By Proposition \ref{prop3.3} this means that $K$, the number of solutions to $a_1 + a_2 = a_3$ with $a_i \in A \cap (H+ x_i)$, satisfies
\[ |K - \alpha_1 \alpha_2 \alpha_3 |H|^2| \leq \epsilon |H|^2,\]
which means that $K \geq \epsilon |H|^2/2$. Thus certainly the number of Schur triples in $A$ is certainly at least $\epsilon |H|^2/2$, which is at least $\epsilon 2^{-2 M(\epsilon) - 1} N^2$. \vspace{11pt}

\noindent The proof of Theorem \ref{reg-thm} is very much in the spirit of the iterative method. One again takes $\mathbf{Struct}$ to be the collection of all subspaces $H \leq \mathbb{F}_2^n$, but here there is no $\mathbf{Config}$. Let $A \subseteq \mathbb{F}_2^n$ be a set, and let $H \in \mathbf{Struct}$. We define the $L^2$-density of $A$ with respect to $H$ by
\[ \delta_H(A) = \mathbb{E}(\frac{|A \cap (H + x)|^2}{|H|^2} \; | \; x \in \mathbb{F}_2^n).\]
In \cite{green-reg} this is called the \textit{index}, and is written $\ind(A;H)$.\vspace{11pt}

\noindent The key to the proof is the following lemma (Lemma 2.2 of \cite{green-reg}), which can be proved by elaborating somewhat on the proof of Lemma \ref{uniformity-hyperplanes} (1).
\begin{lemma}\label{lem2.2} Let $\epsilon \in (0,\frac{1}{2})$ and
suppose that $H \leq \mathbb{F}_2^n$ is a subgroup which is not $\epsilon$-regular for $A$. Then there is a subgroup $H' \leq H$ such that $\mbox{\emph{codim}}(H') \leq 2^{\mbox{\emph{\scriptsize codim}}(H)}$ and $\delta_{H'}(A) \geq \delta_H(A) + \epsilon^3$.
\end{lemma}

\noindent Theorem \ref{reg-thm} is simply a matter of applying Lemma \ref{lem2.2} iteratively. Since $\delta_H(A) \leq 1$ for any $H$, the number of iterations is no more than $\lceil 1/\epsilon^3 \rceil$.\vspace{11pt}

\noindent An unfortunate feature of Theorem \ref{reg-thm} and its proof is that $M(\epsilon)$ grows like a tower of twos of height $\epsilon^{-3}$. This is because each application of Lemma \ref{lem2.2} results in an exponentiation of the codimension of $H$. By adapting a brilliant construction of Gowers \cite{gowers-tower}, which shows that Szemer\'edi's regularity lemma for graphs must have tower type bounds, we were able to show that $M(\epsilon)$ \textit{must} be at least as bad as a tower of twos of height about $\log(1/\epsilon)$.\vspace{11pt}

\noindent We were not able to produce a similar example in the setting of Theorem \ref{almost-sum-free}. 

\begin{problem}
Find a ``reasonable'' bound for $C_2(\delta)$, the quantity appearing in Theorem \ref{almost-sum-free}, or prove that no such bound exists.
\end{problem}

\noindent In fact for $G = \mathbb{F}_2^n$ I am not able to exclude the possibility that $C_2(\delta)$ can be a polynomial in $\delta$. This need not be the case for $G = \mathbb{Z}/N\mathbb{Z}$, due to the Behrend example of a large subset of $\{1,\dots, N\}$ containing no 3-term AP. See \cite{green-reg} for a further discussion.\vspace{11pt}

\noindent The corresponding for graphs (relating to Theorem \ref{ruz-szem-thm}) is also wide open, though again it is known that $C_1(\delta)$ cannot be taken to be polynomial in $\delta$.

\section{From Finite Fields to $\{1,\dots,N\}$}\label{sec8}

\noindent We have now seen several examples concerning additive combinatorics in finite fields. However, for many of the problems we have considered it is an analogue in $\{1,\dots,N\}$ or (more-or-less equivalently) in $\mathbb{Z}/N\mathbb{Z}$ which is actually of interest. \vspace{11pt}

\noindent In recent years the passage from finite fields to the integers, at least for problems concerning configurations of the type we have been discussing in the last four sections, has started to form into something resembling a theory. This is thanks to the work of Bourgain \cite{bourgain2} on finding good bounds for $r_3(N)$. \vspace{11pt}

\noindent Bourgain's ideas are developed in detail in his original paper, of course, and have also been discussed in \cite{green-reg} and \cite{tao-notes}. In this section we restrict ourselves to a few remarks which illustrate the important points.\vspace{11pt}

\noindent Consider the problem of finding a bound for $r_3(G)$ using the iteration method, where $G$ is an abelian group with order $N$ and no 2-torsion. It is not hard to see (essentially by changing the letter $\xi$ to $\gamma$ in Proposition \ref{prop3.3}) that if $A \subseteq G$ has density $\alpha$ and substantially fewer than $\alpha^3 N^2$ 3-term APs then $A$ has a non-trivial large Fourier coefficient, that is to say 
\[ \widehat{A}(\gamma) := \sum_x A(x) \gamma(x)\]
has magnitude a large fraction of $N$ for some non-trivial character $\gamma \in G^{*}$. \vspace{11pt}

\noindent It is not immediately clear how to use this information. We can no longer assert that $A$ has increased density on a subspace, because in a general group $G$ there is no such thing as a subspace. What we can show, rather painlessly, is that $A$ has increased density on a translate of a \textit{Bohr set}, that is to say a set of the form $x + B(\{\gamma\},\epsilon)$, where
\[ B(\{\gamma\},\epsilon) := \{x \in G: |1 - \gamma(x)| \leq \epsilon\}.\]
The reader who has followed the various iterative arguments in the last four sections might now suggest that we define \textbf{Struct} to be the collection of all Bohr sets $B(\Gamma,\epsilon)$, where $\Gamma = \{\gamma_1,\dots,\gamma_d\}$ is a set of characters and 
\[ B(\Gamma,\epsilon) := \{x \in G: |1 - \gamma_j(x)| \leq \epsilon \qquad \mbox{for all $j = 1,\dots,d$}\}.\]
Note that in $\mathbb{F}_3^n$ a Bohr set is the same thing as a subspace when $\epsilon < 1/4$.
Such a strategy is clearly not going to be without its difficulties. If $S \in \mathbf{Struct}$, it looks as though we are going to have to make some sense of what it means to do Fourier analysis on $S$. Since $B(\Gamma,\epsilon)$ is not a group, this will certainly not be a trivial matter. \vspace{11pt}

\noindent In fact, $B = B(\Gamma,\epsilon)$ is quite a long way from being a group. The homomorphism $(\gamma_1,\dots,\gamma_d) : G \rightarrow \mathbb{T}^d$ carries $B$ into a small $d$-dimensional box $D$. If one picks $x,x'$ at random in $D$, the chance that $x + x' \in D$ is just $2^{-d}$. Hence one expects that typically $|B + B| \approx 2^{d}|B|$, which compares unfavourably with the result $|H + H| = |H|$ which holds if $H \leq G$ is a genuine subspace.\vspace{11pt}

\noindent We will not, in this survey, go into the details of what we mean by Fourier analysis on $B$, nor how the large doubling constant of $B$ is unpleasant in this context. We hope the reader will believe us when we say that reducing the doubling constant is a very helpful thing to do.\vspace{11pt}

\noindent Bourgain's advance is to consider $B$ not by itself, but together with another Bohr set $B' := B(\Gamma,\epsilon')$, where $\epsilon'$ is much smaller than $\epsilon$ . Then if $x \in B$ and $x' \in B'$ we have $x + x' \in B(\Gamma,\epsilon + \epsilon')$, a set which ought not to be much larger than $B$. Thus $|B + B'| \approx |B|$, and we may think of the pair $(B,B')$ as behaving like an approximate group. Roughly speaking, it turns out to indeed be possible to run an iterative argument in which \textbf{Struct} is the collection of all such pairs $(B,B')$.\vspace{11pt}

\noindent There are a number of further technicalities to be overcome. One interesting one is that our assertion that $B(\Gamma,\epsilon + \epsilon')$ is not much larger than $B(\Gamma,\epsilon)$ is not true in general. Suppose, for example, that $G = \mathbb{F}_5^n$, that the characters in $\Gamma$ are linearly independent and that $\epsilon < 2\sin (\pi/5), \epsilon + \epsilon' > 2\sin(\pi/5)$. Then $|B(\Gamma,\epsilon)| = 5^{n - d}$, whilst $|B(\Gamma,\epsilon + \epsilon')| = 3^{-d}5^n$. Bourgain circumvents this difficulty by using an averaging argument to show that for a typical $\epsilon$ the size of $B(\Gamma,\epsilon)$ is roughly invariant under small perturbations of $\epsilon$. Tao \cite{tao-notes} observed that one could also replace Bohr sets by \textit{smoothed} Bohr sets, and then such difficulties go away. I implemented this idea slightly differently in \cite{green-reg}, defining the a smoothed Bohr ``set'' by
\[ \widetilde{B}(\Gamma,\epsilon)(x) := \int^{\infty}_0 B(K,t)(x) \frac{e^{-t/\epsilon}}{\epsilon} \, dt.\]

\noindent We conclude this section by giving an up-to-date summary of the extent to which the problems of the last four sections have been given Bourgain's treatment. Of course, in the original paper \cite{bourgain2} the question of $r_3(G)$ was treated (actually, Bourgain only treats $r_3(N)$ but it is clear that his methods work in an arbitrary $G$). In \cite{green-reg} the results of \S \ref{sec7} are all fully generalised to any finite abelian $G$, and in particular Theorem \ref{almost-sum-free} is proved in this general setting. As regards adapting the methods of \S \ref{sec5} to obtain a bound of the form $r_{\angle}(G) \ll N(\log \log N)^{-c}$, this ought to be possible (Shkredov, work in progress). Finally there is the issue of transferring the arguments of \S \ref{sec6} to obtain a bound of the form $r_4(G) \ll N(\log N)^{-c}$. In particular one would like this for $G = \mathbb{Z}/N\mathbb{Z}$, which would imply that $r_4(N) \ll N(\log N)^{-c}$. Since the argument for $r_4(\mathbb{F}_5^n)$ is already rather difficult, one should not expect this to be at all straightforward. Even describing the correct generalisation of the notion of quadratic form to an arbitrary $G$ is not straightforward \cite{green-tao-manifesto}.

\section{Progressions in Sumsets}\label{sec9}

\noindent As promised, we now move onto questions of a somewhat more miscellaneous nature. This section concerns the following problem.
\begin{problem}\label{sumset-aps}
Let $A \subseteq \{1,\dots,N\}$ be a set of size $N/10$ \textup{(}say\textup{)}. Must $A + A$ contain a long arithmetic progression?
\end{problem}
\noindent Bourgain \cite{bourgain1} proved that the answer is ``yes''; $A + A$ must contain a surprisingly long arithmetic progression. If $L(N,\alpha)$ is the smallest $l$ for which there is a set $A \subseteq \{1,\dots,N\}$ of cardinality $\alpha N$ such that $A + A$ does not contain a progression of length $l$, then Bourgain showed that $L(N,1/10) \gg \exp(c(\log N)^{1/3})$. In \cite{green-aps-sumsets} this was improved to $L(N,1/10) \gg \exp(c(\log N)^{1/2})$. An example of Ruzsa \cite{ruzsa-niveau} shows that $L(N,1/10) \ll \exp(c_{\epsilon}(\log N)^{2/3 + \epsilon})$.\vspace{11pt}

\noindent It seems as though the natural finite field analogue of Problem \ref{sumset-aps} involves replacing ``arithmetic progression'' by ``coset of a subspace''.
\begin{problem}\label{sumset-hyperplane}
Write $D(n,\alpha)$ for the smallest $d$ for which there is $A \subseteq \mathbb{F}_2^n$ of density $\alpha$ such that $A + A$ does not contain a coset of a subspace of dimension $d$. Estimate $D(n,\alpha)$.
\end{problem}

\noindent The techniques of \cite{green-aps-sumsets} adapt to this situation in a straightforward manner, and one obtains the following.

\begin{theorem}\label{thm9.3} Suppose that $\alpha \geq n^{-1/4}$. Then $D(n,\alpha) \geq \alpha^2 n/80$.
\end{theorem}

\noindent A detailed proof of this fact may be found in \cite{green-rkpnotes}. In keeping with the philosophy of this survey, some of the details are rather cleaner than in the orginal argument \cite{green-aps-sumsets} which applied to subsets of $\{1,\dots,N\}$.\vspace{11pt}

\noindent A more dramatic difference between the finite field case and the original setting of Problem \ref{sumset-aps} can be observed when one tries to adapt Ruzsa's construction to the finite field setting. 

\begin{theorem}[Ruzsa's niveau sets in $\mathbb{F}_2^n$]\label{ruz-niv}
$D(n,1/4) \leq n - \sqrt{n}$.
\end{theorem} 
\proof 
 Let $A$ be the set of all vectors $x \in \mathbb{F}_2^n$ with at least $n/2 + \sqrt{n}/2$ ones with respect to the standard basis. By the central limit theorem the number of ones in a random vector $(x_1,\dots,x_n)$ is roughly normally distributed with mean $n/2$ and standard deviation $\sqrt{n}/2$, and so for large $n$ we have $|A| \geq 2^{n-2}$. Now any vector $x \in A + A$ must have at least $\sqrt{n}$ zeros. Using this fact, we shall prove that $A + A$ meets all translates of all $(n - \lfloor\sqrt{n}\rfloor)$-dimensional subspaces. Indeed, write $d = \lfloor\sqrt{n}\rfloor$ and suppose that $U$ is a translate of some subspace of dimension $n = d$. $U$ can be written as
\[ U \; = \; \left\{a_0 + \lambda_1a_1 + \dots + \lambda_{n-d}a_{n-d} \; : \; \lambda_i \in \mathbb{F}_2\right\},\] where the $a_i$ are linearly independent. Write $a_i$ in component form as $(a_i^{(j)})_{j = 1}^n$. The column rank of the matrix $(a_{ij})$ is $n - d$, and hence so is the row rank. Without loss of generality, suppose that the first $n - d$ rows $(a_1^{(j)},\dots,a_{n-d}^{(j)})$, $j = 1,\dots,n-d$, are linearly independent. Then we can solve the $n - d$ equations
\[ a_0^{(j)} + \lambda_1a_1^{(j)} + \dots + \lambda_{n-d}a_{n-d}^{(j)} \; = \; 1\] for the $\lambda_i$, giving a vector in $U$ with no more than $d$ zeros.\endproof

\begin{problem}
Narrow the gap between Theorems \ref{thm9.3} and \ref{ruz-niv}.
\end{problem}
\noindent My suspicion is that the upper bound of Theorem \ref{ruz-niv} is closer to the truth.\vspace{11pt}

\noindent I cannot resist mentioning two problems which were raised at the AIM conference on additive combinatorics. The first is due to Croot:

\begin{problem} Fix $\theta \in (0,1)$. What is 
\[ l(\theta) = \limsup_{N \rightarrow \infty} \min_{A \subseteq [N], |A| = N^{1 - \theta}} (\mbox{length of the longest progression in $A + A$})?\]
\end{problem}

\noindent In words, we are interesting in finding subsets $A \subseteq [N]$ with density $N^{-\theta}$ such that $A + A$ contains no long arithmetic progression. Croot states that the bounds $2/\theta - 1 \leq l(\theta) \ll 2^{1/\theta}$ are known. The upper bound comes by considering a multidimensional progression of dimension about $\theta \log_2 N$: it would be interesting to see whether a construction related to niveau sets gives anything better.\vspace{11pt}

\noindent The second question is due to Katznelson:

\begin{problem} What is the measure of the largest open subset $A$ of the torus $\mathbb{T}^d$ for which $A - A$ does not contain a 1-dimensional subgroup? In particular, is it $2^{-d}$?
\end{problem}

\section{Freiman's Theorem}\label{sec10}

\noindent A great deal of the material in this section was communicated to me in person by Imre Ruzsa, and is reproduced here and in the supplementary document \cite{sup-doc-2} (which contains proofs) with his kind permission. The reader will also wish to consult Ruzsa's own survey article \cite{ruzsa-survey}, as well as the material from the AIM conference on Additive Combinatorics \cite{aim-notes}.\vspace{11pt}

\noindent This section concerns Problem \ref{prob2} of the introduction. Let $A \subseteq \mathbb{F}_2^{\infty}$ have doubling at most $K$, meaning that we have the inequality $|A + A| \leq K|A|$. What can be said about the structure of $A$?\vspace{11pt}

\noindent It is hard to think of any examples of sets $A$ with this property other than cosets of subspaces, and large subsets of them. In fact, these are the only such examples as was shown by Imre Ruzsa \cite{ruzsa-freimangroups}. The best known bounds for a result of this type are due to Ruzsa and the author \cite{green-ruzsa3}:

\begin{theorem}[Freiman's theorem in $\mathbb{F}_2^{\infty}$]\label{freiman-torsion} Let $A \subseteq \mathbb{F}_2^{\infty}$ be a finite set with $|A + A| 
\leq K|A|$. Then $A$ is contained within a coset of some subgroup $H \leq \mathbb{F}_2^{\infty}$
with $|H| \leq K^22^{2K^2 -2}|A|$.
\end{theorem}
\noindent A version of this result, with somewhat weaker bounds, will be a consequence of Proposition \ref{ruzsa-prop} below (which is also due to Imre Ruzsa). \vspace{11pt}

\noindent Theorem \ref{freiman-torsion} gives, in a weak sense, a complete description of sets with small doubling. We showed that if $|A + A| \leq K|A|$ then $A$ is contained in a coset of a subspace of size at most $K^22^{2K^2-2}|A|$; conversely, if $A$ has this property then it is clear that $|A + A| \leq K^22^{2K^2 - 2}|A|$. It would be of great interest to have a structure theorem which does not result in exponential losses in $K$ of this sort. Perhaps one can even arrange things so that one has a result of the form

\[ \mbox{doubling constant $K$} \Longrightarrow \mbox{structure} \Longrightarrow \mbox{doubling constant $K'$}, \] where $K'$ is \textit{polynomial} in $K$.\vspace{11pt}

\noindent It is easy to see that such a structure theorem would have to take a form somewhat different from Theorem \ref{freiman-torsion}. Indeed if one takes $A \subseteq \mathbb{F}_2^{\infty}$ to be a subspace $H$ together with $K$ points $x_1,\dots,x_K$ such that $\mbox{Span}(x_1,\dots,x_K) \cap H = \{0\}$ then it is clear that $|A + A| \leq K|A|$, but that the smallest coset-of-a-subspace containing $A$ has size roughly $2^K|A|$.\vspace{11pt}

\noindent Ruzsa \cite{ruzsa-freimangroups} reports that Katalin Marton has suggested that one should be looking for a covering of $A$ by a small number $C_1(K)$ of cosets of some rather smaller subspace of size $C_2(K)|A|$. I agree with this, and it is to some extent believeable that $C_1(K)$ and $C_2(K)$ can be polynomial in $K$. Ruzsa was probably the first to actually dare to conjecture this, and he certainly states such a conjecture explicitly in \cite{ruzsa-survey}. Such matters are also touched upon (in the $\mathbb{Z}$-setting) in \cite{bourgain3,gowers-rough-struct}.\vspace{11pt}

\noindent Imre Ruzsa indicated to me a large part of the following proposition giving a number of statements equivalent to such a structure theorem. The proof may be found in \cite{sup-doc-2}.

\begin{proposition}[Ruzsa]\label{ruzsa-prop}
The following five statements are equivalent.
\begin{enumerate}
\item If $A \subseteq \mathbb{F}_2^{\infty}$ has $|A + A| \leq K|A|$, then there is $A' \subseteq A$, $|A'| \geq |A|/C_1(K)$, which is contained in a coset of some subspace of size at most $C_2(K)|A|$.
\item If $A \subseteq \mathbb{F}_2^{\infty}$ has $|A + A| \leq K|A|$, then $A$ may be covered by at most $C_3(K)$ cosets of some subspace of size at most $C_4(K)|A|$.
\item If $A \subseteq \mathbb{F}_2^{\infty}$ has $|A + A| \leq K|A|$, and if additionally there is a set $B$, $|B| \leq K$, such that $A + B = A + A$, then $A$ may be covered by at most $C_5(K)$ cosets of some subspace of size at most $C_6(K)|A|$.
\item Suppose that $f : \mathbb{F}_2^m \rightarrow \mathbb{F}_2^{\infty}$ is a function with the property that $|\{f(x) + f(y) - f(x+y) : x,y \in \mathbb{F}_2^m\}| \leq K$. Then $f$ may be written as $g + h$, where $g$ is linear and $|\mbox{Im}(h)| \leq C_7(K)$.
\item Suppose that $f : \mathbb{F}_2^m \rightarrow \mathbb{F}_2^{\infty}$ is a function with the property that for at least $2^{3m}/K$ of the quadruples $(x_1,x_2,x_3,x_4) \in \mathbb{F}_2^m$ with $x_1 + x_2 = x_3 + x_4$ we have $f(x_1) + f(x_2) = f(x_3) + f(x_4)$. Then there is an affine linear function $g : \mathbb{F}_2^m \rightarrow \mathbb{F}_2^{\infty}$ such that $f(x) = g(x)$ for at least $2^m/C_8(K)$ values of $x$.
\end{enumerate}
Furthermore if $C_i(K)$ is bounded by a polynomial in $K$ for all $i \in I$, where $I$ is any of the sets $\{1,2\},\{3,4\},\{5,6\},\{7\},\{8\}$ then in fact $C_i(K)$ is bounded by a polynomial in $K$ for all $i$.
\end{proposition}
\noindent\textit{Remarks.} Statement (4) is perhaps the most elegant and natural one here. Observe also that (4) is rather easy with the bound $C_7(K) = 2^K$. Thus Proposition \ref{ruzsa-prop} implies a weak version of Theorem \ref{freiman-torsion}. It is the possibility of polynomial bounds for $C_i(K)$ that is the most interesting feature of this proposition. Let us call this the PFR conjecture:

\begin{conjecture}[Polynomial Freiman-Ruzsa conjecture for $\mathbb{F}_2^n$]
The function $C_7(K)$ \textup{(}and hence all of the other functions $C_i(K)$, $i = 1,\dots,8$\textup{)}, can be taken to be polynomial in $K$.
\end{conjecture}

\noindent The following question has implications for PFR.

\begin{question}
Let $A \subseteq \mathbb{F}_2^n$ be a set of density $\alpha$. Then $2A - 2A$ contains a subspace with codimension $f(\alpha)$. What is the behaviour of $f(\alpha)$?\end{question}

\noindent Using a Fourier-analytic technique of Bogolyubov \cite{bogolyubov} one may show that $f(\alpha) \ll \alpha^{-2}$, and a refinement of this technique due to Chang \cite{chang} allows one to improve this to $f(\alpha) \ll \alpha^{-1}\log (1/\alpha)$. We have not been able to rule out the possibility that $f(\alpha) \ll \log(1/\alpha)$, which if true would imply PFR.\vspace{11pt}

\noindent The proof of Proposition \ref{ruzsa-prop} uses an important result known as Pl\"unnecke's inequality \cite{plunnecke}, a new proof of which was found by Ruzsa \cite{ruzsa-plun}. This states that if $A$ is a subset of any abelian group $G$, and if $|A + A| \leq K|A|$, then we have the inequality $|sA - tA| \leq K^{s+t}|A|$ for any positive integers $s,t$. The reader may observe that (1) of Proposition \ref{ruzsa-prop} implies a much stronger bound for some large subset $A' \subseteq A$, for large $s,t$, at least if there is a good bound on $C_2(K)$. We may call such an $A'$ \textit{subpl\"unnecke}. Nets Katz asked me to formulate a principle to the effect that $A$ being subpl\"unnecke implies that $A$ is very economically contained in some coset of a subspace. The following result is my best effort so far in this direction:

\begin{proposition}\label{subplunprop} Let $A \subseteq \mathbb{F}_2^{\infty}$, and suppose that there is a constant $B$ such that $|tA| \leq t^B|A|$ for all $t \geq B\log B$. Then $A$ is contained in a union of $2^{CB\log B}$ cosets of some subspace having size at most $|A|$.
\end{proposition}

\noindent The hope, of course, is that one might be able to show that if $|A + A| \leq K|A|$ then $A$ has a large subset $A'$ which is subpl\"unnecke in the sense of Proposition \ref{subplunprop}, for some reasonably small $B$ (ideally, $B = O(\log K/\log \log K)$, which would imply PFR).\vspace{11pt}

\noindent For me the most important reason for wanting to understand the PFR conjecture is the implications it would have for our understanding of quadratic Fourier coefficients. In particular, PFR in $\mathbb{F}_5^n$ (the formulation is obvious) would imply a positive solution to the PGI Conjecture (Conjecture \ref{quad-inverse}).

\begin{proposition}\label{quad-prop} Suppose that \textup{PFR} is true in $\mathbb{F}_5^n$. Then \textup{PGI} is true. That is, let $f : \mathbb{F}_5^n \rightarrow [-1,1]$ be a function with $\mathbb{E}f = 0$, and suppose that $\Vert f \Vert_{U^3} \geq \delta$. Then there is a quadratic form $q$ on $\mathbb{F}_5^n$ such that
\[ |\mathbb{E} f(x) \omega^{q(x)}| \gg \delta^{C},\] for some absolute constant $C$.
\end{proposition}

\noindent The deduction is given in \cite{green-tao-manifesto}.\vspace{11pt}

\noindent In my opinion it would be very interesting to determine whether PGI has any implications for PFR. It is just plausible that this represents the most natural way to attack PFR, though at the moment we have little idea how to carry out such a programme.\vspace{11pt}

\noindent The results of this section may be discussed in the context of general abelian groups $G$. However, the issues are of a rather different nature to those discussed in \S \ref{sec8}. Freiman's original work concerned subsets of $\mathbb{Z}$, and was quite geometric in feel. See \cite{bilu,freiman,green-edinburgh-notes} for a further discussion. Ruzsa's proof \cite{ruzsa-gaps} has proved much more adaptable, and recently Ruzsa and the author \cite{green-ruzsa4} were able to obtain a structure theorem for sets with small doubling which is valid in any abelian group.

\begin{theorem}[G. -- Ruzsa]
Let $G$ be an abelian group, and suppose that $A \subseteq G$ has $|A + A| \leq K|A|$. The $A$ is contained in a set of the form $H + P$, where $H$ is a subgroup, $P$ is a generalised arithmetic progression, the dimension of $P$ is $\leq C_9(K)$ and $|H||P| \leq C_{10}(K)$.
\end{theorem}
\noindent\textit{Remark.} A generalised arithmetic progression of dimension $d$ is a set of the form
\[ \{ a_0 + \lambda_1 a_1 + \dots + \lambda_d a_d : 0 \leq \lambda_i \leq L_i \qquad \mbox{for $i = 1,\dots,d$}.\}\]
We obtain the bounds $C_9(K) \ll K^{C}$ and $C_{10}(K) \ll e^{K^{C}}$, for some absolute constant $C$.

\providecommand{\bysame}{\leavevmode\hbox to3em{\hrulefill}\thinspace}
\providecommand{\MR}{\relax\ifhmode\unskip\space\fi MR }
\providecommand{\MRhref}[2]{%
  \href{http://www.ams.org/mathscinet-getitem?mr=#1}{#2}
}
\providecommand{\href}[2]{#2}

     \end{document}